\def\version{0.18}

\def\journal{XXX}
\def\titlep{Non-existence of universal $R$-matrix
for some C$^{*}$-bialgebras}
\documentclass[11pt]{article}
\usepackage{graphicx,ifthen}
\usepackage{amssymb}
\usepackage{amsmath}

\font\germ=eufm10 at12pt

\def\goth#1{\hbox{\germ#1}}

\newcommand{\qed}{\hbox{\rule[-2pt]{3pt}{6pt}}}
\newcommand{\qedh}{\hfill\qed \\}

\newcommand{\vv}{\vspace{.3in}}

\newtheorem{Thm}{Theorem}[section]

\newtheorem{rem}[Thm]{Remark}

\newtheorem{defi}[Thm]{Definition}
\newtheorem{lem}[Thm]{Lemma}

\newtheorem{prob}[Thm]{Problem}

\newcommand{\kn}{\Large\bf
$K\hspace{-.4cm} N$
\Large\bf\vv }

\def\cal#1{\mathcal #1}
\def\con{{\cal O}_{n}}

\def\edot{=1,\ldots,n}
\def\pr{{\it Proof.}\quad}

\def\coa{{\cal O}_{A}}
\def\co#1{{\cal O}_{#1}}
\def\disp#1{{\displaystyle #1}}
\def\ck{Cuntz-Krieger}

\def\brl{branching law}
\def\bfsnl{{\rm BFS}_{N}(\Lambda)}

\addtocounter{footnote}{1}
\def\sftt#1{
\setcounter{equation}{0}
\addtocounter{footnote}{1}
\section{#1}
}

\def\ssft#1{\subsection{#1}}

\begin{document}
%
%
\def\autherp{Katsunori Kawamura}
\def\emailp{e-mail: kawamura@kurims.kyoto-u.ac.jp.}
\def\addressp{{\small {\it College of Science and Engineering,
 Ritsumeikan University,}}\\
{\small {\it 1-1-1 Noji Higashi, Kusatsu, Shiga 525-8577, Japan}}
}

\def\infw{\Lambda^{\frac{\infty}{2}}V}
\def\zhalfs{{\bf Z}+\frac{1}{2}}
\def\ems{\emptyset}
\def\pmvac{|{\rm vac}\!\!>\!\! _{\pm}}
\def\vac{|{\rm vac}\rangle _{+}}
\def\dvac{|{\rm vac}\rangle _{-}}
\def\ovac{|0\rangle}
\def\tovac{|\tilde{0}\rangle}
\def\expt#1{\langle #1\rangle}
\def\zph{{\bf Z}_{+/2}}
\def\zmh{{\bf Z}_{-/2}}
\def\brl{branching law}
\def\bfsnl{{\rm BFS}_{N}(\Lambda)}
\def\scm#1{S({\bf C}^{N})^{\otimes #1}}
\def\mqb{\{(M_{i},q_{i},B_{i})\}_{i=1}^{N}}
\def\zhalf{\mbox{${\bf Z}+\frac{1}{2}$}}
\def\zmha{\mbox{${\bf Z}_{\leq 0}-\frac{1}{2}$}}
\newcommand{\mline}{\noindent
\thicklines
\setlength{\unitlength}{.1mm}
\begin{picture}(1000,5)
\put(0,0){\line(1,0){1250}}
\end{picture}
\par
 }
\def\ptimes{\otimes_{\varphi}}
\def\delp{\Delta_{\varphi}}
\def\delps{\Delta_{\varphi^{*}}}
\def\gamp{\Gamma_{\varphi}}
\def\gamps{\Gamma_{\varphi^{*}}}
\def\sem{{\sf M}}
\def\hdelp{\hat{\Delta}_{\varphi}}
\def\tilco#1{\tilde{\co{#1}}}
\def\ndm#1{{\bf M}_{#1}(\{0,1\})}
\def\cdm#1{{\cal M}_{#1}(\{0,1\})}
\def\tndm#1{\tilde{{\bf M}}_{#1}(\{0,1\})}
\def\sck{{\sf CK}_{*}}
\def\hdel{\hat{\Delta}}
\def\ba{\mbox{\boldmath$a$}}
\def\bb{\mbox{\boldmath$b$}}
\def\bc{\mbox{\boldmath$c$}}
\def\bd{\mbox{\boldmath$d$}}
\def\be{\mbox{\boldmath$e$}}
\def\bk{\mbox{\boldmath$k$}}
\def\bm{\mbox{\boldmath$m$}}
\def\bp{\mbox{\boldmath$p$}}
\def\bq{\mbox{\boldmath$q$}}
\def\bu{\mbox{\boldmath$u$}}
\def\bv{\mbox{\boldmath$v$}}
\def\bw{\mbox{\boldmath$w$}}
\def\bx{\mbox{\boldmath$x$}}
\def\by{\mbox{\boldmath$y$}}
\def\bz{\mbox{\boldmath$z$}}
\def\bomega{\mbox{\boldmath$\omega$}}
\def\N{{\bf N}}
\def\lxm{L_{2}(X,\mu)}
\def\qtimes{\otimes_{\tilde{\varphi}}}
\def\ul#1{\underline{#1}}
\def\titlepage{

\noindent
{\bf 
\noindent
\thicklines
\setlength{\unitlength}{.1mm}
\begin{picture}(1000,0)(0,-300)
\put(0,0){\kn \knn\, for \journal\, Ver.\version}
\put(0,-50){\today}
\end{picture}
}
\vspace{-2.3cm}
\quad\\
{\small file: \textsf{tit01.txt,\, J1.tex}
\footnote{
${\displaystyle
\mbox{directory: \textsf{\fileplace}, 
file: \textsf{\incfile},\, from \startdate}}$}}
\quad\\
\framebox{
\begin{tabular}{ll}
\textsf{Title:} &
\begin{minipage}[t]{4in}
\titlep
\end{minipage}
\\
\textsf{Author:} &\autherp
\end{tabular}
}
{\footnotesize	
\tableofcontents }
}

%
%
%
\setcounter{section}{0}
\setcounter{footnote}{0}
\setcounter{page}{1}
\pagestyle{plain}

%
%
\title{\titlep}
\author{\autherp\thanks{\emailp}
\\
\addressp}
\date{}
\maketitle

%
%
\begin{abstract}
For a C$^{*}$-bialgebra $A$ with a comultiplication $\Delta$,
a universal $R$-matrix of $(A,\Delta)$ is defined as 
a unitary element in the multiplier algebra $M(A\otimes A)$
of $A\otimes A$ which is an intertwiner between 
$\Delta$ and its opposite comultiplication $\Delta^{op}$.
We show that there exists no universal $R$-matrix
for some C$^{*}$-bialgebras.
\end{abstract}

\noindent
{\bf Mathematics Subject Classifications (2000).} 16W30, 81R50. \\
\\
{\bf Key words.} universal $R$-matrix, C$^{*}$-bialgebra

%
%
\sftt{Introduction}
\label{section:first}
We have studied C$^{*}$-bialgebras and their representations.
In this paper, we consider a universal $R$-matrix of a C$^{*}$-bialgebra.
Since a C$^{*}$-bialgebra is not always a bialgebra
in a sense of the purely algebraic theory of quantum groups,
we generalize the definition of universal $R$-matrix to C$^{*}$-bialgebras.
Next, we consider whether a certain C$^{*}$-bialgebra has a universal
$R$-matrix or not.
In this section, we show our motivation,
definitions and main theorem.

%
%
\ssft{Motivation}
\label{subsection:firstone}
In this subsection, we roughly explain our motivation 
and the background of this study.
Explicit mathematical definitions will 
be shown after $\S$ \ref{subsection:firsttwo}.

Let $\co{*}$ denote 
the direct sum of all Cuntz algebras except $\co{\infty}$:
%
%
\begin{equation}
\label{eqn:cuntztwo}
\co{*}=\co{1}\oplus\co{2}\oplus\co{3}\oplus\co{4}\oplus\cdots
\end{equation}
where $\co{1}$ denotes the $1$-dimensional C$^{*}$-algebra ${\bf C}$
for convenience.
In \cite{TS02}, 
we constructed a non-cocommutative comultiplication $\Delta_{\varphi}$ 
of ${\cal O}_{*}$.
The C$^{*}$-bialgebra $(\co{*},\delp)$ has no antipode (with a dense domain).
With respect to $\delp$,
tensor product formulae of representations of $\con$'s
are well studied \cite{TS01,TS07}.
Details about $(\co{*},\delp)$ are 
will be explained in $\S$ \ref{subsection:firstthree}.

On the other hand, 
in the theory of quantum groups,
a universal $R$-matrix for a quasi-cocommutative bialgebra 
is important for an application to mathematical physics 
and low-dimensional topology \cite{Drinfeld,Jimbo,Kassel}.
Therefore, it is meaningful for a given bialgebra
to find its universal $R$-matrix if it exists.
However, C$^{*}$-bialgebra is not always a bialgebra in
a sense of the theory of purely algebraic case \cite{EAbe, Kassel}.
Therefore,
notions in purely algebraic case without change
can not be always applied to C$^{*}$-bialgebra (see also Remark \ref{rem:weak}(ii)).
Related studies were also considered by Van Daele and Van Keer 
for Hopf $*$-algebras \cite{VanDaeleVanKeer}.

Our interests are to define a notion 
of universal $R$-matrix of a C$^{*}$-bialgebra and  
to clarify whether 
$(\co{*},\delp)$ in (\ref{eqn:cuntztwo}) has a universal $R$-matrix
or not.
In this paper, we show the negative result, that is,
$(\co{*},\delp)$  has no universal $R$-matrix.
For this purpose,
we show a statement about the non-existence of universal $R$-matrix of
a general C$^{*}$-bialgebra (Lemma \ref{lem:maintwob}).
 
%
%
\ssft{Definitions}
\label{subsection:firsttwo}
In this subsection,
we recall definitions of C$^{*}$-bialgebra,
and introduce a universal $R$-matrix of a C$^{*}$-bialgebra.
For two C$^{*}$-algebras $A$ and $B$,
let ${\rm Hom}(A,B)$ and $A\otimes B$ 
denote the set of all $*$-homomorphisms from $A$ to $B$ and 
the minimal C$^{*}$-tensor product of $A$ and $B$, respectively.
Let $M(A)$ denote the multiplier algebra of $A$.

At first,
we prepare terminologies about C$^{*}$-bialgebra according 
to \cite{TS02,KV,MNW}.
A pair $(A,\Delta)$ is a {\it C$^{*}$-bialgebra}
if $A$ is a C$^{*}$-algebra and $\Delta\in {\rm Hom}(A,M(A\otimes A))$ 
such that $\Delta$ is nondegenerate and the following holds:
%
%
\begin{equation}
\label{eqn:bialgebratwo}
(\Delta\otimes id)\circ \Delta=(id\otimes\Delta)\circ \Delta.
\end{equation}
We call $\Delta$ the {\it comultiplication} of $A$.
Remark that we do not assume $\Delta(A)\subset A\otimes A$.
Furthermore,
$A$ has no unit in general for a C$^{*}$-bialgebra $(A,\Delta)$.
From these,
a C$^{*}$-bialgebra is not always a bialgebra in
a sense of the purely algebra theory \cite{EAbe,Kassel}.

According to \cite{Drinfeld,Kassel,VanDaeleVanKeer},
we introduce a unitary universal $R$-matrix and the quasi-cocommutativity
for a C$^{*}$-bialgebra as follows.
%
%
\begin{defi}
\label{defi:rmat}
Let $(A,\Delta)$ be a  C$^{*}$-bialgebra.
\begin{enumerate}
\item
The map $\tilde{\tau}_{A,A}$ from $M(A\otimes A)$ to $M(A\otimes A)$
defined as
%
%
\begin{equation}
\label{eqn:tilde}
\tilde{\tau}_{A,A}(X)(x\otimes y)\equiv \tau_{A,A}(X(y\otimes x))
\quad(X\in M(A\otimes A),\,x\otimes y\in A\otimes A)
\end{equation}
is called the extended flip
where $\tau_{A,A}$ denotes the flip of $A\otimes A$.
\item
The map $\Delta^{op}$ from $A$ to $M(A\otimes A)$
defined as
%
%
\begin{equation}
\label{eqn:opposite}
\Delta^{op}(x)\equiv \tilde{\tau}_{A,A}(\Delta(x))\quad(x\in A)
\end{equation}
is called the opposite comultiplication of $\Delta$.
%
\item
A C$^{*}$-bialgebra $(A,\Delta)$ is cocommutative
if $\Delta=\Delta^{op}$. 
\item
An element $R\in M(A\otimes A)$
is called a (unitary) universal $R$-matrix of $(A,\Delta)$
if $R$ is a unitary
and 
%
%
\begin{equation}
\label{eqn:univ}
R\Delta(x)R^{*}=\Delta^{op}(x)\quad
(x\in A).
\end{equation}
In this case,
we state that $(A,\Delta)$ is quasi-cocommutative.
\end{enumerate}
\end{defi}
%
%
\begin{rem}
\label{rem:weak}
{\rm 
\begin{enumerate}
\item
The additional assumption of unitarity of 
a universal $R$-matrix is natural
for $*$-algebras.
\item
If $A$ is unital,
then $M(A\otimes A)=A\otimes A$ and
$\tilde{\tau}_{A,A}=\tau_{A,A}$.
In addition, if $(A,\Delta)$ is quasi-cocommutative with
a universal $R$-matrix $R$, then $R\in A\otimes A$.
In the purely algebraic theory of quantum groups,
a bialgebra has a unit by definition \cite{Kassel}.
Hence the quasi-cocommutativity makes sense 
by using a universal $R$-matrix as an invertible element 
in the tensor square of the bialgebra.
On the other hand, 
there is no unit in C$^{*}$-bialgebra in general by definition.
If there is no unit, then
there is no invertible element in the algebra.
Hence the quasi-cocommutativity and universal $R$-matrix make no sense 
for C$^{*}$-bialgebras if one uses the purely algebraic axiom without change.
\item
If $(A,\Delta)$ is cocommutative, then
the unit of $M(A\otimes A)$ is a universal $R$-matrix of $(A,\Delta)$.
Hence $(A,\Delta)$ is quasi-cocommutative.
On the other hand,
if a quasi-cocommutative C$^{*}$-bialgebra $(A,\Delta,R)$ is not 
cocommutative,
then $R$ is not a scalar multiple of the unit of $M(A\otimes A)$.
Therefore the non-quasi-commutativity 
is stronger than the non-commutativity
as same as the purely algebraic case.
\item
About examples of universal $R$-matrix in purely algebraic theory,
see examples in Chap. VIII.2 of \cite{Kassel}.
\end{enumerate}
}
\end{rem}

%
%
\ssft{Main theorem}
\label{subsection:firstthree}
In this subsection, we show our main theorem.
Before that,
we explain the C$^{*}$-bialgebra
$(\co{*},\delp)$ in \cite{TS02}  more closely.
Let $\con$ denote the Cuntz algebra for $2\leq n<\infty$ \cite{C},
that is, the C$^{*}$-algebra which is universally generated by
generators $s_{1},\ldots,s_{n}$ satisfying
$s_{i}^{*}s_{j}=\delta_{ij}I$ for $i,j\edot$ and
$\sum_{i=1}^{n}s_{i}s_{i}^{*}=I$
where $I$ denotes the unit of $\con$.
The Cuntz algebra $\con$ is simple, that is,
there is no nontrivial two-sided closed ideal.
This implies that any unital representation of $\con$ is faithful.

Redefine the C$^{*}$-algebra $\co{*}$
as the direct sum of the set $\{\con:n\in {\bf N}\}$ of Cuntz algebras:
%
%
\begin{equation}
\label{eqn:cuntbi}
\co{*}\equiv \bigoplus_{n\in {\bf N}} \con
=\{(x_{n}):\|(x_{n})\|\to 0\mbox{ as }n\to\infty\}
\end{equation}
where ${\bf N}=\{1,2,3,\ldots\}$
and $\co{1}$ denotes the $1$-dimensional C$^{*}$-algebra for convenience.
For $n\in {\bf N}$,
let $s_{1}^{(n)},\ldots,s_{n}^{(n)}$ denote
canonical generators of $\con$
where $s_{1}^{(1)}$ denotes the unit of $\co{*}$.
For $n,m\in {\bf N}$,
define the embedding $\varphi_{n,m}$ of $\co{nm}$
into $\con\otimes \co{m}$ by
%
%
\begin{equation}
\label{eqn:embeddingone}
\varphi_{n,m}(s_{m(i-1)+j}^{(nm)})\equiv s_{i}^{(n)}\otimes s_{j}^{(m)}
\quad(i=1,\ldots,n,\,j=1,\ldots,m).
\end{equation}
For the set $\varphi\equiv \{\varphi_{n,m}:n,m\in {\bf N}\}$ in
(\ref{eqn:embeddingone}),
define the $*$-homomorphism $\delp$ from $\co{*}$ to $\co{*}\otimes \co{*}$ by
%
%
\begin{eqnarray}
\label{eqn:dpone}
\delp\equiv& \oplus\{\delp^{(n)}:n\in {\bf N}\},\\
\nonumber
\\
\label{eqn:dptwo}
\delp^{(n)}(x)\equiv &\disp{\sum_{(m,l)\in {\bf N}^{2},\,ml=n}\varphi_{m,l}(x)}
\quad(x\in \con,\,n\in {\bf N}).
\end{eqnarray}
Then  the following holds:
The pair $(\co{*},\delp)$ is a non-cocommutative C$^{*}$-bialgebra
(\cite{TS02}, Theorem 1.1);
There is no antipode for any dense subbialgebra of $(\co{*},\delp)$
(\cite{TS02}, Theorem 1.2(v)).
About much further properties of $(\co{*},\delp)$, see \cite{TS02,TS04}.
About a generalization of $(\co{*},\delp)$, see \cite{TS05}.

Our main theorem is stated as follows.
%
%
\begin{Thm}
\label{Thm:noexist}
There is no universal $R$-matrix of $(\co{*},\delp)$.
\end{Thm}

In \cite{Drinfeld},
Drinfel'd constructed a universal $R$-matrix
by taking a completion of a given bialgebra
with respect to a certain topology
(see also Chap. XVI of \cite{Kassel}).
As an analogy of this,
we propose the following problem
for non-quasi-cocommutative C$^{*}$-bialgebras.
%
%
\begin{prob}
\label{prob:one}
{\rm
For a non-quasi-cocommutative C$^{*}$-bialgebra $(A,\Delta)$
(for example, $(\co{*},\delp)$),
construct an extension $(\tilde{A},\tilde{\Delta})$
of $(A,\Delta)$ such that 
$(\tilde{A},\tilde{\Delta})$ is a quasi-cocommutative.
}
\end{prob}

\noindent
As a related topic, 
Drinfel'd's quantum double \cite{Drinfeld}
is known as a method  of construction of a braided Hopf algebra
from any finite dimensional Hopf algebra with invertible antipode.
Difficulties of application of this method to infinite dimensional case
are discussed in \cite{VanDaeleVanKeer}. 

In $\S$ \ref{section:second},
we will prove Theorem \ref{Thm:noexist}.
In $\S$ \ref{section:third}, we will show another example 
of non-quasi-cocommutative C$^{*}$-bialgebra.
%
%
\sftt{Proof of Theorem \ref{Thm:noexist}}
\label{section:second}
We prove Theorem \ref{Thm:noexist} in this section.
For this purpose,
we recall multiplier algebra and nondegenerate homomorphism,
and show a lemma.

%
%
\ssft{Multiplier algebra and nondegenerate homomorphism}
\label{subsection:secondone}
For a C$^{*}$-algebra $A$,
let $A^{''}$ denote the enveloping von Neumann algebra of $A$.
The {\it multiplier algebra} $M(A)$ of $A$ is defined by
%
%
\begin{equation}
\label{eqn:ma}
M(A)\equiv 
\{a\in A^{''}:aA\subset A,\,
Aa\subset A\}.
\end{equation}
Then $M(A)$ is a unital C$^{*}$-subalgebra of $A^{''}$.
Especially,
$A=M(A)$ if and only if $A$ is unital.
The algebra $M(A)$ is the completion of $A$
with respect to the strict topology.

A $*$-homomorphism from $A$ to $B$ is not always extended to the map 
from $M(A)$ to $M(B)$.
If $({\cal H},\pi)$ is a nondegenerate representation of $A$,
that is,
$\pi(A){\cal H}$ is dense in ${\cal H}$,
then there exists a unique extension of $\pi$ in ${\rm Hom}(M(A),{\cal B}({\cal H}))$.
We state that $f\in {\rm Hom}(A,M(B))$ is {\it nondegenerate} if $f (A)B$ 
is dense in a C$^{*}$-algebra $B$.
If both $A$ and $B$ are unital and $f$ is unital,
then $f$ is nondegenerate. 
For a $*$-homomorphism $f$ from $A$ to $M(B)$,
if $f$ is nondegenerate, then 
$f$ is called a {\it morphism} from $A$ to $B$ \cite{W3}.
If $f\in{\rm Hom}(A,B)$ is nondegenerate,
then we can regard $f$ as a morphism from $A$ to $B$
by using the canonical embedding of $B$ into $M(B)$.
Each morphism $f$ from $A$ to $B$ can be extended uniquely 
to a homomorphism $\tilde{f}$ from $M(A)$ to $M(B)$
such that $\tilde{f}(m) f(b)a = f(mb)a$ 
for $m\in M(B), b\in B$, and $a \in A$.
If $f$ is injective, then so is $\tilde{f}$.

%
%
\ssft{A lemma}
\label{subsection:secondtwo}
For a unitary element $U$ in a unital $*$-algebra ${\goth A}$,
we define the (inner) $*$-automorphism ${\rm Ad}U$ of ${\goth A}$  by
${\rm Ad}U(x)\equiv UxU^{*}$ for $x\in {\goth A}$.
%
%
%
\begin{lem}
\label{lem:maintwob}
Let $(A,\Delta)$ be a C$^*$-bialgebra.
If $(A,\Delta)$ is quasi-cocommutative,
then 
for any two nondegenerate representations
$\pi_{1}$ and $\pi_{2}$ of $A$,
$(\pi_{1}\otimes \pi_{2})\circ \Delta$ 
and
$(\pi_{2}\otimes \pi_{1})\circ \Delta$ 
are unitarily equivalent
where we write the extension of
$\pi_{i}\otimes \pi_{j}$ on $M(A\otimes A)$
as $\pi_{i}\otimes \pi_{j}$ for $i,j=1,2$.
\end{lem}
%
%
\pr
Since $\pi_{ij}\equiv 
\pi_{i}\otimes \pi_{j}$ is also nondegenerate,
we can extend $\pi_{ij}$ on $M(A\otimes A)$
and use the same symbol $\pi_{ij}$ for its extension
for $i,j=1,2$.
From this, $\pi_{ij}\circ \Delta$ is well-defined.
Let ${\cal H}_{i}$ denote the representation space of $\pi_{i}$ for $i=1,2$
and 
let $T$ denote the flip between
${\cal H}_{1}\otimes {\cal H}_{2}$
and ${\cal H}_{2}\otimes {\cal H}_{1}$,
which is a unitary operator.
Then the following holds:
%
%
\begin{equation}
\label{eqn:assumebb}
\pi_{21}\circ \Delta
= {\rm Ad}T\circ  \pi_{12}\circ \Delta^{op}.
\end{equation}

By assumption,
there exists a universal $R$-matrix $R$ of $(A,\Delta)$. 
Define $S\equiv \pi_{12}(R)$. 
From (\ref{eqn:assumebb}),
%
%
\begin{equation}
\label{eqn:uniunib}
\pi_{21}\circ \Delta
=
{\rm Ad}T\circ\pi_{12}\circ ({\rm Ad}R\circ \Delta)
=
{\rm Ad}W\circ \pi_{12}\circ  \Delta
\end{equation}
where $W\equiv TS$.
This means that 
$(\pi_{2}\otimes \pi_{1})\circ \Delta$ and 
$(\pi_{1}\otimes \pi_{2})\circ \Delta$
are unitarily equivalent.
\qedh

\noindent
From Lemma \ref{lem:maintwob},
we see that
the study of tensor products of representations of a bialgebra $(A,\Delta)$
is useful to verify whether $(A,\Delta)$ is quasi-cocommutative or not.
The idea of Lemma \ref{lem:maintwob}
is a modification of a well-known fact in the purely algebraic case.
We explain it as follows:
Let $A$ be a bialgebra with a comultiplication $\Delta$
in a sense of purely algebraic theory \cite{Kassel}.
If $(A,\Delta)$ has a universal $R$-matrix 
(where the assumption of ``braided" is not necessary), then 
for any two representations $\pi_{1},\pi_{2}$ of $A$,
$(\pi_{1}\otimes \pi_{2})\circ \Delta$ and 
$(\pi_{2}\otimes \pi_{1})\circ \Delta$ are equivalent
(\cite{Kassel}, Proposition VIII.3.1(a)).
In order to define 
$(\pi_{1}\otimes \pi_{2})\circ \Delta$ and 
$(\pi_{2}\otimes \pi_{1})\circ \Delta$ for the C$^{*}$-bialgebraic case,
we assume that two representations are nondegenerate
because
neither $\pi_{1}\otimes \pi_{2}$ 
nor $\pi_{2}\otimes \pi_{1}$ 
can be always extended on $M(A\otimes A)$ without any assumption.

%
%
\ssft{Proof of Theorem \ref{Thm:noexist}}
\label{subsection:secondthree}
Let $P_{n}$ denote the canonical projection of $\co{*}$
to $\con$.
We identify a representation
$\pi$ of $\con$ 
with the representation $\pi\circ P_{n}$ of $\co{*}$ here.
If $\pi$ is a nondegenerate representation of $\con$,
then 
$\pi$ is also a nondegenerate representation of $\co{*}$.
If two representations $\pi_{1}$ and $\pi_{2}$ of $\con$ are  
not unitarily equivalent, 
then $\pi_{1}$ and $\pi_{2}$ are also not as 
representations of $\co{*}$.

It is known that
there exist two unital representations $\pi_{1}$ and $\pi_{2}$ of $\co{2}$
(especially, $\pi_{1}$ and $\pi_{2}$ are nondegenerate)
such that $(\pi_{1}\otimes \pi_{2})\circ \delp|_{\co{4}}$
and $(\pi_{2}\otimes \pi_{1})\circ \delp|_{\co{4}}$
are not unitarily equivalent from Example 4.1 of \cite{TS01}.
Hence $(\pi_{1}\otimes \pi_{2})\circ \delp$
and $(\pi_{2}\otimes \pi_{1})\circ \delp$
are not unitarily equivalent.
From these and the contraposition of Lemma \ref{lem:maintwob},
the statement holds.
\qedh

%
%
\sftt{Example}
\label{section:third}
We show another example of non-quasi-cocommutative C$^{*}$-bialgebra.
We call a matrix $A$ {\it nondegenerate} 
if any column and any row are not zero.
For $1\leq n<\infty$,
let $\ndm{n}$ denote the set of all
nondegenerate $n\times n$ matrices with entries $0$ or $1$.
In particular, $\ndm{1}=\{1\}$.
Define 
%
%
\begin{equation}
\label{eqn:nondegenerate}
\ndm{*}\equiv \cup \{\ndm{n}:n\geq 1\}.
\end{equation}
For $A=(a_{ij})\in \ndm{n}$,
let 
$\coa$ denote the \ck\ algebra associated with $A$ \cite{CK}.
Define the direct sum ${\sf CK}_{*}$ of $\{\coa:A\in \ndm{*}\}$:
%
%
\begin{equation}
\label{eqn:ckatwo}
{\sf CK}_{*}\equiv \bigoplus_{A\in \ndm{*}}\coa
\end{equation}
where we define $\co{1}\equiv {\bf C}$ for convenience.
Then ${\sf CK}_{*}$ is a C$^{*}$-bialgebra (\cite{TS05}, Theorem 1.2)
such that $\co{*}$ in (\ref{eqn:cuntztwo})
is a C$^{*}$-subbialgebra of ${\sf CK}_{*}$
and it is also a direct sum component of ${\sf CK}_{*}$ 
(\cite{TS05}, Theorem 1.3(iv)).
From this,
Lemma \ref{lem:maintwob} and Theorem \ref{Thm:noexist},
the C$^{*}$-bialgebra ${\sf CK}_{*}$ is not quasi-cocommutative.


%
%

\end{document}